\documentclass[10pt]{article}

\usepackage{amsmath}
\usepackage{amsfonts}
\usepackage{amssymb}
\usepackage[english]{babel}
\usepackage{amsthm}
\usepackage{extarrows}
\usepackage{amscd}
\usepackage{tikz}
\usepackage{mathdots}
\usetikzlibrary{arrows,shapes}

\def\pf{\par\noindent {\bf Proof}~\par\noindent}

\newcommand{\mR}{\mathbb{R}}

\newcommand{\mN}{\mathbb{N}}

\newcommand{\mE}{\mathbb{E}}
\newcommand{\mS}{\mathbb{S}}

\newcommand{\mT}{\mathbb{T}}
\newcommand{\mU}{\mathbb{U}}

\newcommand{\mcH}{\mathcal{H}}

\newcommand{\mcD}{\mathcal{D}}

\newcommand{\mcV}{\mathcal{V}}
\newcommand{\mcW}{\mathcal{W}}

\newcommand{\ux}{\underline{x}}

\newcommand{\uom}{\underline{\omega}}

\newcommand{\p}{\partial}
\newcommand{\dirac}{\underline{\p}}

\newcommand{\eop}{\hfill$\square$}

\newcommand{\onehalf}{\frac{1}{2}}

\newcommand{\invr}{\frac{1}{r}}
\newcommand{\invrsq}{\frac{1}{r^2}}

\newtheorem{lemma}{Lemma}
\newtheorem{proposition}{Proposition}
\newtheorem{definition}{Definition}
\newtheorem{remark}{Remark}
\newtheorem{corollary}{Corollary}
\newtheorem{property}{Property}
\newtheorem{example}{Example}

\begin{document}

\title{On the radial derivative of the delta distribution}
\author{Fred Brackx$^\ast$, Frank Sommen$^\ast$ \& Jasson Vindas$^\ddagger$}

\date{\small{$^\ast$ Department of Mathematical Analysis, 
Faculty of Engineering and Architecture, Ghent University\\
$^\ddagger$ Department of Mathematics, Faculty of Sciences, Ghent University}}

\maketitle

\begin{center}
{\em Dedicated to our co-author Frank on the occasion of his 60th birthday.}
\end{center}
\vspace{5mm}


\begin{abstract}
Possibilities for defining the radial derivative of the delta distribution $\delta(\ux)$ in the setting of spherical coordinates are explored.  This leads to the introduction of a new class of continuous linear functionals similar to but different from the standard distributions. The radial derivative of $\delta(\ux)$ then belongs to that new class of so-called signumdistributions. It is shown that these signumdistributions obey easy-to-handle calculus rules which are in accordance with those for the standard distributions in $\mR^m$.
\end{abstract}


\vspace{5mm}


\section{Introduction}


Given a distribution $T$ in $\mR^m$ expressed in spherical coordinates $\ux = r \uom, r =|\ux|, \uom \in \mS^{m-1}$, $\mS^{m-1}$ being the unit sphere in $\mR^m$, one may wonder if a meaning could be given to the actions $\p_r \, T$, $r\, T$, $\uom\, T$ and the like, which, in principle, are `forbidden actions" and thus {\em a priori} not defined in the standard setting. Indeed, differentiation of distributions is well--defined with respect to the standard cartesian coordinates $x_1, x_2, \ldots, x_m$ and multiplication of a distribution is only allowed by a smooth function, a condition which, clearly, is not satisfied by the functions $r$ and $\uom$. We will tackle this problem for one specific distribution: the delta distribution $\delta(\ux)$, and in the first place concentrate on a possible definition of its radial derivative $\p_r\, \delta(\ux)$.\\

\noindent
The delta distribution is pointly supported at the origin, it is rotation invariant: $$\delta(A\, \ux) = \delta(\ux), \ \forall \, A \in {\rm SO}(m)$$ it is even:  $\delta(-\ux) = \delta(\ux)$ and it is homogeneous of order $(-m)$: $$\delta(a \ux) = \frac{1}{|a|^m}\, \delta(\ux)$$ So in a first, naive, approach, one could think of $\p_r\, \delta(\ux)$ as a distribution which remains pointly supported at the origin, rotation invariant, even and homogeneous of degree $(-m-1)$.  Temporarily leaving aside the even character, on the basis of the other cited characteristics the distribution $\p_r\, \delta(\ux)$ should take the following form:
$$
\p_r\, \delta(\ux) = c_0\, \p_{x_1} \delta(\ux) + \cdots + c_m\, \p_{x_m} \delta(\ux)
$$
and it becomes immediately clear that this approach is impossible since all distributions appearing in this decomposition are odd and not rotation invariant, whereas $\p_r\, \delta(\ux)$ is assumed to be even and rotation invariant. It could be that $\p_r\, \delta(\ux)$ is either the zero distribution or is no longer pointly supported at the origin, both possibilities being quite unacceptable.  But another idea is that $\p_r\, \delta(\ux)$ is not a usual distribution anymore, and the same for $r\, \delta(\ux)$ and $\uom\, \delta(\ux)$.  In Section \ref{alternative} we will introduce a new class of bounded linear functionals on an appropriate space of test functions, very similar to but different from the standard distributions, and we will show that the above cited three distributions belong to that new class.


\section{The delta distribution in cartesian coordinates}


The delta distribution $\delta(\ux) \in \mcD'(\mR^m)$ is, quite naturally, very well--known and frequently used in physics to model point sources in various field theories. Let us summarize its properties.  It is a scalar distribution defined by $$\langle \  \delta(\ux) , \varphi(\ux) \  \rangle = \varphi(0), \ \forall \varphi \in \mcD(\mR^m)$$ which is of finite order zero, with cartesian derivatives given by
$$\langle \  \p_{x_j}^s \delta(\ux) , \varphi(\ux) \ \rangle = (-1)^s \langle \   \delta(\ux) , \p_{x_j}^s \varphi(\ux) \ \rangle =  (-1)^s \{\p_{x_j}^s \varphi(\ux)\}|_{\ux=0}, \ j = 1,\ldots,m$$ In particular the action of the Dirac operator $\dirac = \sum_{j=1}^m \, e_j\, \p_{x_j}$ results into the vector--valued ditribution given by
$$
\langle \ \dirac\, \delta(\ux)\, , \varphi(\ux) \ \rangle = - \, \langle \  \delta(\ux)\, , \dirac\, \varphi(\ux) \ \rangle = - \, \{\dirac \varphi(\ux)\}|_{\ux=0}
$$
Note that we are using here the basis vectors $(e_j, j=1,\ldots,m)$ of $\mR^m$ as Clifford vectors, generating the Clifford algebra $\mR_{0,m}$, for which $e_j^2 = - 1, e_i \wedge e_j = e_ie_j = - e_je_i = - e_j \wedge e_i, e_i \cdot e_j = 0, i \neq j=1,\ldots,m$. For more on Clifford algebras we refer to e.g.\ \cite{port}. In this way the Dirac operator, which may be seen as a Stein--Weiss projection of the gradient operator (see e.g.\ \cite{stein}) and underlies the higher dimensional theory of monogenic functions (see e.g.\ \cite{dss}),  linearizes the Laplace operator: $\dirac^2 = -\, \Delta$, which is the Fischer or Fourier dual to $\ux^2 = -\, |\ux|^2$, $\ux$ being the Clifford vector variable $\ux = \sum_{j=1}^m \, e_j\, x_j$.\\
The action of the Euler operator $$\mE = -\, \ux \cdot \dirac = \sum_{j=1}^m \, x_j\, \p_{x_j}$$ on the delta distribution reveals the latter's homogeneous character: $\mE \, \delta(\ux) = (-m)\, \delta(\ux)$, while the action of the (bivector--valued) angular momentum operator $$\Gamma = -\, \ux \wedge \dirac = -\, \sum_{j<k}\, e_j e_k\, (x_j \, \p_{x_k} - x_k \, \p_{x_j} )$$ leads to $\Gamma\, \delta(\ux) = 0$.


\section{The delta distribution in spherical coordinates}


Introducing spherical coordinates: $\ux = r \, \uom, r= |\ux|, \uom = \sum_{j=1}^m\, e_j \omega_j \in \mS^{m-1}$, the Dirac operator takes the form 
$$\dirac = \uom \, \p_r +  \invr \, \p_{\uom}$$ To give an idea how the angular differential operator $\p_{\uom} = \sum_{j=1}^m\, e_j \p_{\omega_j}$ looks like, we mention here its explicit form in dimension $m=2: \p_{\uom} = e_{\theta} \, \p_{\theta}$ and in dimension $m=3: \p_{\uom} = e_{\theta} \, \p_{\theta} + e_{\varphi} \, \frac{1}{\sin{\theta}}\, \p_{\varphi}$, where the meaning of the angular coordinates $\theta$ and $\varphi$ is straightforward.
The Euler operator in spherical coordinates reads: $\mE = r \, \p_r$, while the angular momentum operator $\Gamma$ takes the form $\Gamma = - \uom \, \p_{\uom} = -\, \uom \wedge \p_{\uom}$. In Clifford analysis (see e.g.\ \cite{dss}) this operator $\Gamma$ is mostly called the spherical Dirac operator and the fact that $\Gamma\, \delta(\ux) = 0$ confirms that the delta distribution $\delta(\ux)$, when expressed in spherical coordinates, only depends on the radial distance $r$, in other words the delta distribution is spherically symmetric. Finally the Laplace operator is written in spherical coordinates as
$$\Delta = \p_r^2 + (m-1)\,  \invr \, \p_r + \invrsq\, \Delta^*$$ where $\Delta^* = \uom \p_{\uom} - \p_{\uom} \p_{\uom} = -\, \p_{\uom} \cdot \p_{\uom}$ is the Laplace--Beltrami operator, containing only angular derivatives, which implies that $\Delta^* \, \delta(\ux) = 0$ and so $$\Delta \, \delta(\ux) = \p_r^2  \, \delta(\ux) + (m-1)\,  \invr \, \p_r  \, \delta(\ux)$$\\
Using the basic formulae
$$
\mE\, \delta(\ux) = (r\, \p_r)\, \delta(\ux) = (-m)\, \delta(\ux) \quad \mbox{and} \quad \p_{\uom}\, \delta(\ux) = 0
$$
we will now establish the formulae concerning the ``legal'' actions on the delta distribution of two specific differential operators containing the radial derivative, viz. $(\p_r^2)$ and $\left(\invr \, \p_r\right)$. Note that the latter operator is, up to a constant factor, nothing else but the derivative with respect to $r^2$: $\left(\invr \, \p_r\right) = 2\, \p_{r^2}$

\begin{proposition}
\label{partslaplace}
The actions of the operators $(\p_r^2)$ and $\left(\invr \, \p_r\right)$ on the delta distribution are well--defined and it holds that
$$
(\p_r^2)\, \delta(\ux) = \onehalf (m+1) \, \Delta \, \delta(\ux) \quad \mbox{and} \quad \left(\invr \, \p_r\right) \delta(\ux) = -\, \onehalf \, \Delta \, \delta(\ux) 
$$
\end{proposition}

\pf
Applying twice the Euler operator we obtain consecutively
\begin{eqnarray*}
(r \, \p_r)^2 \, \delta & = & m^2 \, \delta\\
(r \, \p_r) \, \delta + r^2 \, (\p_r^2) \, \delta & = & m^2 \, \delta\\
r^2 \, (\p_r^2) \, \delta & = & m(m+1) \, \delta\\
(\p_r^2) \, \delta & = & m(m+1) \, \invrsq \delta\\
(\p_r^2) \, \delta & = & \onehalf (m+1) \, \Delta \, \delta  + S_1
\end{eqnarray*}
with $S_1 = c_0 \, \delta + \sum_{j=1}^m \, c_j\, \p_{x_j} \, \delta$, since $r^2 \, \Delta \, \delta = (x_1^2 + \cdots + x_m^2) (\p_{x_1}^2 + \cdots + \p_{x_m}^2) \, \delta = 2 m \, \delta$, $r^2 \, \delta = 0$ and $r^2 \, \p_{x_j} \delta = 0$.\\
On the other hand, we have
\begin{eqnarray*}
\left(\invr \, \p_r\right)(r \, \p_r) \, \delta & = & (-m) \, \left(\invr \, \p_r\right) \delta\\
\left(\invr \, \p_r\right) \delta +  (\p_r^2) \, \delta & = & (-m) \, (\invr \, \p_r) \, \delta\\
\left(\invr \, \p_r\right) \delta & = & - \frac{1}{m+1} \, (\p_r^2)\delta\\
\left(\invr \, \p_r\right) \delta & = & - \onehalf \, \Delta \, \delta  - \frac{1}{m+1} \, S_1
\end{eqnarray*}
Now invoking the Laplace operator, we obtain
\begin{eqnarray*}
\Delta \, \delta & = & (\p_r^2) \, \delta + (m-1) \left(\invr \, \p_r\right) \delta\\
                       & = & \left( \onehalf (m+1) \, \Delta \, \delta  + S_1\right) + (m-1) \, \left( - \onehalf \, \Delta \, \delta  - \frac{1}{m+1} \, S_1 \right)\\
                       & =  & \Delta \, \delta + \frac{2}{m+1} \, S_1
\end{eqnarray*}
from which it follows that  $S_1 = 0$.\eop

\begin{remark}
{\rm
Note that the results of Proposition \ref{partslaplace} are consistent with the expression of the Laplace operator in spherical coordinates since
$$
(\p_r^2)  \, \delta(\ux) + (m-1)\left(\invr \, \p_r\right) \delta(\ux) = \onehalf (m+1) \, \Delta \, \delta(\ux) + (m-1) (-\, \onehalf)\, \Delta \, \delta(\ux) = \Delta \, \delta(\ux)
$$
}
\end{remark}

\begin{proposition}
\label{omegadr}
The action of the operator $(\uom\, \p_r)$ on the delta distribution is well--defined and it holds that
$$
(\uom\, \p_r)\, \delta(\ux) = \dirac\, \delta(\ux) = \sum_{j=1}^m\, e_j \, \p_{x_j}  \delta(\ux)
$$
\end{proposition}

\pf
The result easily follows from the spherical form of the Dirac operator, taking into account that $\p_{\uom}\, \delta(\ux) = 0$. \eop

\begin{remark}
{\em
The result of Proposition \ref{omegadr}, which may also been written componentwise as
$$
(\omega_j \p_r)\, \delta(\ux) = \p_{x_j} \delta(\ux), \quad j=1,\ldots,m
$$
is consistent with the expression of the Laplace operator in spherical coordinates since
\begin{eqnarray*}
\Delta\, \delta(\ux) &=& -\, \dirac^2\, \delta(\ux) = -\, (\uom\, \p_r + \invr\, \p_{\uom})(\uom\, \p_r)\, \delta(\ux)\\ 
&=& \p_r^2\, \delta(\ux) -  \p_{\uom}( \uom) \, \invr \p_r \, \delta(\ux)\\
&=& \p_r^2  \, \delta(\ux) + (m-1)\,  \invr \, \p_r  \, \delta(\ux)
\end{eqnarray*}
where the well--known result: $\p_{\uom}(\uom) = 1-m$ was taken into account.
}
\end{remark}

By iterated application of the formulae established in the Propositions \ref{partslaplace} and \ref{omegadr}, the following results are obtained.

\begin{corollary}
\label{omegaformulae}
For all $k \in \mN$ one has
\begin{eqnarray*}
(\uom\, \p_r)^{2k}\, \delta(\ux) & = & (-1)^k\, \p_r^{2k}\, \delta(\ux) = \frac{1}{2^k k!}\, (m+1)(m+3)\cdots(m+2k-1)\, \dirac^{2k}\, \delta(\ux)\\
(\uom\, \p_r)^{2k+1}\, \delta(\ux) & = & (-1)^k\, \uom\, \p_r^{2k+1}\, \delta(\ux) = \frac{1}{2^k k!}\, (m+1)(m+3)\cdots(m+2k-1)\, \dirac^{2k+1}\, \delta(\ux)\\[1mm]
(\uom\, \p_r)\, \dirac^{2k}\, \delta(\ux) & = & \dirac^{2k+1}\, \delta(\ux) = \dirac^{2k}\, (\uom\, \p_r)\, \delta(\ux)\\[1mm]
(\uom\, \p_r)\, \dirac^{2k+1}\, \delta(\ux) & = & \frac{m+2k+1}{2(k+1)}\, \dirac^{2k+2}\, \delta(\ux) = \frac{m+2k+1}{2(k+1)}\,  \dirac^{2k+1}\, (\uom\, \p_r)\, \delta(\ux)
\end{eqnarray*}
\end{corollary}

\begin{corollary}
For all $k \in \mN$ one has
$$
\left(\invr\, \p_r\right)^k \delta(\ux) = \frac{1}{2^k k!}\, \dirac^{2k}\, \delta(\ux) = (-1)^k\, \frac{1}{2^k k!}\, \Delta^k\, \delta(\ux)
$$
\end{corollary}

\noindent
The formulae obtained in Corollary \ref{omegaformulae} may be generalized by considering products of radial derivatives of the delta distribution and natural powers of  the radial distance. The results of the following Proposition \ref{omegaformulaegeneralized} are obtained by a straightforward computation invoking the identities, with $k \geq \ell$,

\begin{eqnarray*}
(i) \qquad \quad \ux^{2\ell}\, \dirac^{2k}\, \delta(\ux) &=& (2k)(2k-2)\cdots(2k-2\ell+2)\\
&&\hspace{23mm}(m+2k-2)(m+2k-4)\cdots(m+2k-2\ell)\, \dirac^{2k-2\ell}\delta(\ux)\\
(ii) \qquad \ux^{2\ell+1}\, \dirac^{2k}\, \delta(\ux) &=& (2k)(2k-2)\cdots(2k-2\ell)\\
&&\hspace{20mm}(m+2k-2)(m+2k-4)\cdots(m+2k-2\ell)\, \dirac^{2k-2\ell-1}\delta(\ux)\\
(iii) \qquad \ux^{2\ell}\, \dirac^{2k+1}\, \delta(\ux) &=& (2k)(2k-2)\cdots(2k-2\ell+2)\\
&&\hspace{20mm}(m+2k)(m+2k-2)\cdots(m+2k-2\ell+2)\, \dirac^{2k-2\ell+1}\delta(\ux)\\
(iv) \ \  \ux^{2\ell+1}\, \dirac^{2k+1}\, \delta(\ux) &=& (2k)(2k-2)\cdots(2k-2\ell+2)\\
&&\hspace{30mm}(m+2k)(m+2k-2)\cdots(m+2k-2\ell)\, \dirac^{2k-2\ell}\delta(\ux)
\end{eqnarray*}

\noindent
It turns out that when the sum of the order of the radial derivative $\p_r$ and the power of $r$ is even then the product is a well--defined scalar radial operator, and when this sum is odd then the product is a well--defined vector operator involving the function $\uom$. 

\begin{proposition}
\label{omegaformulaegeneralized}
One has, with $k \geq \ell$,
\begin{eqnarray*}
(i) \qquad \quad r^{2\ell}\, \p_r^{2k}\, \delta(\ux) &=& (-1)^{k+\ell}\, \frac{1}{2^{k-\ell} (k-\ell)!}\, (m+1)(m+3)\cdots(m+2k-1)\\
&&\hspace{25mm}(m+2k-2)(m+2k-4)\cdots(m+2k-2\ell)\, \dirac^{2k-2\ell}\delta(\ux)\\
(ii) \ \ \ \uom\, r^{2\ell+1}\, \p_r^{2k}\, \delta(\ux) &=& (-1)^{k+\ell}\, \frac{1}{2^{k-\ell-1} (k-\ell-1)!}\, (m+1)(m+3)\cdots(m+2k-1)\\
&&\hspace{22mm}(m+2k-2)(m+2k-4)\cdots(m+2k-2\ell)\, \dirac^{2k-2\ell-1}\delta(\ux)\\
(iii) \ \ \ \uom\, r^{2\ell}\, \p_r^{2k+1}\, \delta(\ux) &=& (-1)^{k+\ell}\, \frac{1}{2^{k-\ell} (k-\ell)!}\, (m+1)(m+3)\cdots(m+2k-1)\\
&&\hspace{22mm}(m+2k)(m+2k-2)\cdots(m+2k-2\ell+2)\, \dirac^{2k-2\ell+1}\delta(\ux)\\
(iv) \quad r^{2\ell+1}\, \p_r^{2k+1}\, \delta(\ux) &=& (-1)^{k+\ell+1}\, \frac{1}{2^{k-\ell} (k-\ell)!}\, (m+1)(m+3)\cdots(m+2k-1)\\
&&\hspace{32mm}(m+2k)(m+2k-2)\cdots(m+2k-2\ell)\, \dirac^{2k-2\ell}\delta(\ux)
\end{eqnarray*}
\end {proposition}


\section{Another attempt to define $\p_r\, \delta(\ux)$}
\label{firstattempt}


First note that the formula
$$
(\uom\, \p_r)\, \delta(\ux) = \dirac\, \delta(\ux) = \sum_{j=1}^m\, e_j \, \p_{x_j}  \delta(\ux)
$$
obtained in Proposition \ref{omegadr} cannot be used to define the radial derivative $\p_r\, \delta(\ux)$ of the delta distribution since multiplication of a distribution by $\uom$ is not allowed.
But the formula
\begin{equation}
\label{drsquared}
(-\, \p_r^2)\, \delta(\ux) = \onehalf (m+1) \, (- \Delta) \, \delta(\ux) 
\end{equation}
obtained in Proposition \ref{partslaplace}, could offer a possibility to define $\p_r\, \delta(\ux)$. Indeed, by taking ``square roots'' we obtain
\begin{equation}
\label{hilbertdr}
H_r \p_r\, \delta(\ux) = \sqrt{\frac{m+1}{2}} \, (- \Delta)^\onehalf \, \delta(\ux) 
\end{equation}
Let us explain this result in more detail. At the right hand--side of (\ref{hilbertdr}) appears the so--called {\em square root of the negative Laplace operator} $(-\Delta)^\onehalf$ which is the convolution operator given, for an appropriate function or distribution $F(\ux)$, by
$$
(-\Delta)^\onehalf [F(\ux)] = -\, \frac{2}{a_{m+1}}\, \textup{Fp}\frac{1}{r^{m+1}} * F(\ux)
$$
the convolution kernel $\textup{Fp}\frac{1}{r^{m+1}}$ being a ``Finite Part" distribution in $\mR^m$, and $a_{m+1} = \frac{2 \pi^{\frac{m+1}{2}}}{\Gamma(\frac{m+1}{2})}$ being the area of the unit sphere $\mS^m$  in $\mR^{m+1}$. Repeated action of the operator $(-\Delta)^\onehalf$ results into the well--known result
$$
(-\Delta)^\onehalf\, (-\Delta)^\onehalf [F] = \left(\frac{2}{a_{m+1}}\right)^2\, \textup{Fp}\frac{1}{r^{m+1}} * [\textup{Fp}\frac{1}{r^{m+1}} * F] = (-\Delta)[F]
$$
recovering in this way for $F(\ux) = \delta(\ux)$, up to the constant $\onehalf (m+1)$, the right hand--side of (\ref{drsquared}). Note that the right hand--side of (\ref{hilbertdr}) reduces to
$$
\sqrt{\frac{m+1}{2}} \, (- \Delta)^\onehalf \, \delta(\ux) =  -\, \frac{2}{a_{m+1}}\,  \sqrt{\frac{m+1}{2}}\, \textup{Fp}\frac{1}{r^{m+1}} * \delta = -\, \sqrt{\frac{m+1}{2}}\, \frac{\Gamma(\frac{m+1}{2})}{\pi^{\frac{m+1}{2}}}\,\textup{Fp}\frac{1}{r^{m+1}}
$$
At the left hand--side of (\ref{hilbertdr}) appears the one--dimensional Hilbert transform $H_r$ in the real variable $r$, given, for an appropriate function or distribution $f(r)$, by
$$
H_r[f(r)] = \frac{1}{\pi}\,  \textup{Pv}\frac{1}{r} * f
$$
the Hilbert kernel $\textup{Pv}\frac{1}{r}$ being a "Principal Value" distribution in $\mR$. As is well--known the one--dimensional Hilbert transform is a linear endomorphism both of the space $\mcD_{L_p}$ of test functions given by
$$
\mcD_{L_p} = \{\phi(t) \in C^{\infty}(\mR): \phi^{(k)}(t) \in L_p(\mR), \forall k \in \mN  \}, \quad 1 < p < +\infty
$$
and of its dual $\mcD'_{L_p}$, with inverse $H^{-1}_r = - H_r$, so that $H_r^2 = - \bf{1}$. Moreover the Hilbert transform is commuting with derivation so that repeated action of $H_r \p_r$ leads to
$$
(H_r \p_r) (H_r \p_r)\, \delta(\ux) = H_r^2 \p_r^2 \, \delta(\ux) = (-\, \p_r^2) \, \delta(\ux)
$$
in this way recovering the left hand--side of (\ref{drsquared}).\\
Now rewriting (\ref{hilbertdr}) as
\begin{equation}
\label{hilbertdr2}
(H_r \p_r)\, \delta(\ux) = \sqrt{\frac{m+1}{2}} \, (- \Delta)^\onehalf \, \delta(\ux) = -\, \sqrt{\frac{m+1}{2}}\, \frac{\Gamma(\frac{m+1}{2})}{\pi^{\frac{m+1}{2}}}\,\textup{Fp}\frac{1}{r^{m+1}}
\end{equation}
it becomes clear that the operator $(H_r\p_r)$ acting on the delta distribution flattens out the point support of $\delta(\ux)$ to the whole of $\mR^m$, and the result 
$(- \Delta)^\onehalf \, \delta(\ux)$ is the only distribution, up to a constant, which is rotation invariant and homogeneous of degree $(-m-1)$. Note by the way that the same phenomenon concerning the support of the delta distribution occurs under the action of the so--called Hilbert--Dirac operator $(H\dirac)$ (see e.g. \cite{hilbertdirac}):
$$
(H\dirac)[\delta(\ux)] =  (- \Delta)^\onehalf \, \delta(\ux) = - \frac{\Gamma(\frac{m+1}{2})}{\pi^{\frac{m+1}{2}}}\,\textup{Fp}\frac{1}{r^{m+1}}
$$
where $H$ stands for the Hilbert transform in $\mR^m$ given, for an appropriate function or distribution $f$, by
$$
H[f] = \mcH * f = -\,  \frac{\Gamma(\frac{m+1}{2})}{\pi^{\frac{m+1}{2}}}\,\textup{Fp}\frac{\uom}{r^{m}} * f
$$
In order to retrieve from (\ref{hilbertdr2}) an expression for $\p_r\, \delta(\ux)$, we can act with the radial Hilbert transform $H_r$ on both sides, leading to
\begin{equation}
\label{dr}
\p_r\, \delta(\ux) =  \sqrt{\frac{m+1}{2}}\, \frac{\Gamma(\frac{m+1}{2})}{\pi^{\frac{m+1}{2}}}\, H_r\left[ \textup{Fp}\frac{1}{r^{m+1}}\right]
\end{equation}
which requires an appropriate definition of the {\em radial Hilbert transform of a distribution} and of a rotation--invariant distribution in particular, a question which is of the same nature as the quest for an acceptable definition of the radial derivative  $\p_r\, \delta(\ux)$. Anyway, as we will show that $\p_r\, \delta(\ux)$ belongs to a new class of so--called signumdistributions, by (\ref{dr}), the radial Hilbert transform of a distribution should belong to that class too.


\section{A physics approach to the delta distribution}


In physics texts one often encounters the following expression for the delta distribution in spherical coordinates:
\begin{equation}
\label{physicsdelta}
\delta(\ux) = \frac{1}{a_m}Ê\, \frac{\delta(r)}{r^{m-1}}
\end{equation}
Apparently this can be explained in the following way. Write the action of the delta distribution as an integral:
\begin{eqnarray*}
\varphi(0) = \langle \ \delta(\ux) , \varphi(\ux) \ \rangle &=& \int_{\mR^m} \, \delta(\ux) \, \varphi(\ux) \, dV(\ux)\\
&=& \int_{0}^{\infty} \, r^{m-1} \delta(\ux) \, dr \, \int_{\mS^{m-1}} \, \varphi(r \, \uom) \, dS_{\uom}\\
&=& a_m\, \int_{0}^\infty \, r^{m-1} \, \delta(\ux) \, \Sigma^{0}[\varphi](r) \, dr
\end{eqnarray*}
using the so--called {\em spherical mean} of the test function $\varphi$ given by
$$
\Sigma^{0}[\varphi](r) = \frac{1}{a_m}Ê\, \int_{\mS^{m-1}} \, \varphi(r \, \uom) \, dS_{\uom}
$$
where $a_m = \frac{2\pi^{\frac{m}{2}}}{\Gamma(\frac{m}{2})}$ is the area of the unit sphere $\mS^{m-1}$ in $\mR^m$.
As it is easily seen that $\Sigma^0[\varphi](0) = \varphi(0)$ it follows that
$$
a_m\, \int_{0}^\infty \, r^{m-1} \, \delta(\ux) \, \Sigma^{0}[\varphi](r) \, dr = \langle \  \delta(r) ,  \Sigma^{0}[\varphi](r) \ \rangle = \int_{0}^\infty \, \delta(r) \, \Sigma^{0}[\varphi](r) \, dr
$$
which explains (\ref{physicsdelta}). However we prefer to interpret this expression as
\begin{equation}
\label{sigmanul}
\varphi(0) = \langle \ \delta(\ux) , \varphi(\ux) \ \rangle = \langle \  \delta(r) ,  \Sigma^{0}[\varphi](r) \ \rangle = \Sigma^{0}[\varphi](0)
\end{equation}
which can be generalized to higher even order Dirac--derivatives of the delta distribution:
$$
\{\dirac^{2\ell}\varphi(\ux)\}|_{\ux=0} = \langle \ \dirac^{2\ell}\,\delta(\ux) , \varphi(\ux) \ \rangle = (-1)^\ell\, C(\ell)\, \langle \  \p_r^{2\ell}\, \delta(r) ,  \Sigma^{0}[\varphi](r) \ \rangle =  (-1)^\ell\, C(\ell)\, \{ \p_r^{2\ell}\,\Sigma^{0}[\varphi](r) \}|_{r=0}
$$
with
$$
C(\ell) = \frac{2^{2\ell}\,\ell!}{(2\ell)!}\left(\frac{m}{2}+\ell-1\right)\cdots\left(\frac{m}{2}\right) =  \frac{2^{2\ell}\,\ell!}{(2\ell)!}\frac{\Gamma(\frac{m}{2}+\ell)}{\Gamma(\frac{m}{2})}, \quad \ell=0,1,2,\ldots
$$
Note that the spherical mean $\Sigma^{0}[\varphi](r)$ is an even function of $r$, whose odd order derivatives vanish at the origin:
$$
 \langle \  -\,\p_r^{2\ell+1}\, \delta(r) ,  \Sigma^{0}[\varphi](r) \ \rangle =  \{\p_r^{2\ell+1}\,\Sigma^{0}[\varphi](r)\}|_{r=0} = 0
$$
For expressing the higher odd order Dirac--derivatives of the delta distribution in a similar way we have to invoke the so--called {\em spherical mean of the second kind} $\Sigma^{1}[\varphi]$, which was introduced in \cite{bds}:
$$
\Sigma^{1}[\varphi](r) = \frac{1}{a_m}Ê\, \int_{\mS^{m-1}} \, \uom\, \varphi(r \, \uom) \, dS_{\uom}
$$
This spherical mean $\Sigma^{1}[\varphi](r)$ is a vector--valued odd function of $r$, whose even order derivatives vanish at the origin:
$$
 \langle \  \p_r^{2\ell}\, \delta(r) ,  \Sigma^{1}[\varphi](r) \ \rangle =  \{\p_r^{2\ell}\,\Sigma^{1}[\varphi](r)\}|_{r=0} = 0
$$
It holds that
\begin{equation}
\label{sigmaeen}
 \langle \ \dirac^{2\ell+1}\,\delta(\ux) , \varphi(\ux) \ \rangle = (-1)^\ell\, C(\ell+1)\, \langle \  \p_r^{2\ell+1}\, \delta(r) ,  \Sigma^{1}[\varphi](r) \ \rangle  
\end{equation}
or
$$
\{\dirac^{2\ell+1}\varphi(\ux)\}|_{\ux=0} = (-1)^\ell\, C(\ell+1)\, \{ \p_r^{2\ell+1}\,\Sigma^{1}[\varphi](r) \}|_{r=0}
$$
In the physics language these results would then be written as
\begin{eqnarray*}
\dirac^{2\ell}\,\delta(\ux) &=& (-1)^\ell\, C(\ell)\, \frac{1}{a_m}\, \frac{\delta^{(2\ell)}(r)}{r^{m-1}}\\
\dirac^{2\ell+1}\,\delta(\ux) &=& (-1)^\ell\, C(\ell+1)\, \frac{1}{a_m}\, \frac{\delta^{(2\ell+1)}(r)}{r^{m-1}}\, \uom
\end{eqnarray*}
which, by means of the results of Corollary \ref{omegaformulae}, lead to the expressions
\begin{eqnarray*}
\p_r^{2\ell}\,\delta(\ux) &=& \frac{1}{(2\ell)!}(m)(m+1)\cdots(m+2\ell-1)\, \frac{1}{a_m}\, \frac{\delta^{(2\ell)}(r)}{r^{m-1}}\\
\uom\, \p_r^{2\ell+1}\,\delta(\ux) &=& \frac{1}{(2\ell+1)!}(m)(m+1)\cdots(m+2\ell)\, \frac{1}{a_m}\, \frac{\delta^{(2\ell+1)}(r)}{r^{m-1}}\, \uom
\end{eqnarray*}
Again $\p_r\, \delta(\ux)$ escapes from this approach.


\section{Spherical representation of a distribution}


When expressing a test function $\varphi(\ux) \in \mcD(\mR^m)$ in spherical coordinates, one obtains a function $\widetilde{\varphi}(r, \uom) = \varphi(r\uom) \in 
\mcD(\mR \times \mS^{m-1})$, but it is clear that not all functions $\widetilde{\varphi}(r, \uom)  \in \mcD(\mR \times \mS^{m-1})$ stem from a test function in $\mcD(\mR^m)$.
However a one--to--one correspondence may be established between the usual space of test functions $\mcD(\mR^m)$ and a specific subspace of $\mcD(\mR \times \mS^{m-1})$.

\begin{lemma} (see \cite{helgason})
\label{iso}
There is a one--to--one correspondence $\varphi(\ux) \leftrightarrow \widetilde{\varphi}(r, \uom) = \varphi(r\uom)$ between the spaces $\mcD(\mR^m)$ and $\mcV = \{\phi(r,\uom) \in \mcD(\mR \times \mS^{m-1})  : \phi$ is even, i.e. $\phi(-r,-\uom) = \phi(r,\uom)$, and  $\{\p_r^n \, \phi(r,\uom) \}|_{r=0}$ is a homogeneous polynomial of degree
$n$ in $(\uom_1,\ldots,\uom_m), \forall n \in \mN\}$.
\end{lemma}

Clearly $\mcV$ is a closed (but not dense) subspace of $\mcD(\mR \times \mS^{m-1})$ and even of $\mcD_E(\mR \times \mS^{m-1})$, where the suffix $E$ refers to the even character of the test functions in that space, and $\mcV$ is endowed with the induced topology of $\mcD(\mR \times \mS^{m-1})$.\\

The one--to--one correspondence between the spaces of test functions $\mcD(\mR^m)$ and $\mcV$ translates into a one--to--one correspondence between the standard distributions $T \in \mcD'(\mR^m)$ and the bounded linear functionals in $\mcV'$; this correspondence is given  by
$$
\langle \ T(\ux) , \varphi(\ux) \ \rangle = \langle \ \widetilde{T}(r,\uom) , \widetilde{\varphi}(r,\uom) \ \rangle
$$

By Hahn--Banach's theorem the bounded linear functional $\widetilde{T}(r,\uom)  \in \mcV'$ may be extended to the distribution $\mT(r,\uom) \in \mcD'(\mR \times \mS^{m-1})$; such an extension is called a {\em spherical representation} of the distribution $T$ (see e.g.\ \cite{vindas}). As the subspace $\mcV$ is not dense in  $\mcD(\mR \times \mS^{m-1})$, the spherical representation of a distribution is {\em not unique}, but if $\mT_1$ and $\mT_2$ are two different spherical representations of the same distribution $T$, their restrictions to $\mcV$ coincide: 
$$
\langle \ \mT_1(r,\uom) , \widetilde{\varphi}(r,\uom) \ \rangle = \langle \ \mT_2(r,\uom) , \widetilde{\varphi}(r,\uom) \ \rangle = 
\langle \ \widetilde{T}(r,\uom) , \varphi(r\uom) \ \rangle = \langle \ T(\ux) , \varphi(\ux) \ \rangle
$$

\noindent
For test functions in $\mcD(\mR \times \mS^{m-1})$ the spherical variables $r$ and $\uom$ are ordinary variables, and thus smooth functions. It follows that for distributions in $\mcD'(\mR \times \mS^{m-1})$ multiplication by $r$ and by $\uom$ and differentiation with respect to $r$ and to $\uom$ are standard well--defined operations, and so
$$
\langle \ \p_r\, \mT(r,\uom) , \Xi(r,\uom) \ \rangle = -\, \langle \ \mT(r,\uom) , \p_r\,  \Xi(r,\uom) \ \rangle 
$$
for all test functions $\Xi(r,\uom) \in \mcD(\mR \times \mS^{m-1})$, and similar expressions for $\p_{\uom}\, \mT$, $r\, \mT$ and $\uom\, \mT$.
However if $\mT_1$ and $\mT_2$ are two different spherical representations of the same distribution $T \in \mcD'(\mR^m)$, then, upon restricting to test functions $\widetilde{\varphi}(r,\uom) \in \mcV$, we are stuck with
$$
-\, \langle \ \mT_1(r,\uom) , \p_r\, \widetilde{\varphi}(r,\uom)  \ \rangle  \neq -\, \langle \ \mT_2(r,\uom) , \p_r\, \widetilde{\varphi}(r,\uom)  \ \rangle 
$$
since $\p_r\, \widetilde{\varphi}(r,\uom)$ is an odd function in the variables $(r,\uom)$ and does no longer belong to $\mcV$ (and neither do $\p_{\uom}\, \widetilde{\varphi}(r,\uom)$, $r\, \widetilde{\varphi}(r,\uom)$ and $\uom\, \widetilde{\varphi}(r,\uom)$). The conclusion is that the concept of spherical representation of a distribution does not allow for an unambiguous definition of the actions proposed, confirming our statement that the solution of our problem lays outside the world of traditional distributions.\\
At the same time it becomes clear why the actions of the operators $\p_r^2$, $\invr\, \p_r$ and $\uom\, \p_r$ on a standard distribution are well-defined instead. Indeed, we have e.g. 
$$
\langle \ \p_r^2\, \mT(r,\uom) , \Xi(r,\uom) \ \rangle =  \langle \ \mT(r,\uom) , \p_r^2\,  \Xi(r,\uom) \ \rangle 
$$
where now $\p_r^2\,  \Xi(r,\uom)$ belongs to $\mcD_E(\mR \times \mS^{m-1})$ which enables restriction to test functions in  $\mcV$ in an unambiguous way.


\section{An alternative class of distributions}
\label{alternative}


As already remarked in the preceding section, $\uom$ is an ordinary (vector) variable in $\mR \times \mS^{m-1}$, whence it makes sense to consider the following subspace of vector--valued test functions in $\mR \times \mS^{m-1}$:
$$
\mcW  = \uom\, \mcV \subset \mcD_O(\mR \times \mS^{m-1}; \mR^m) \subset \mcD(\mR \times \mS^{m-1};\mR^m)
$$
where now the suffix $O$ refers to the odd character of the test functions under consideration, i.e. $\psi(-r,-\uom) = -\, \psi(r,\uom), \forall \psi \in  \mcD_O(\mR \times \mS^{m-1}; \mR^m)$. This space $\mcW$ is endowed with the induced topology of $\mcD(\mR \times \mS^{m-1};\mR^m)$. By definition there is a one--to--one correspondence between the spaces $\mcV$ and $\mcW$.\\

\noindent
For each $\mU(r,\uom) \in \mcD'(\mR \times \mS^{m-1};\mR^m)$ we define $\widetilde{U}(r,\uom) \in \mcW'$ by the restriction
$$
\langle \ \widetilde{U}(r,\uom) , \uom\, \widetilde{\varphi}(r,\uom) \ \rangle =  \langle \ \mU(r,\uom) , \uom\, \widetilde{\varphi}(r,\uom) \ \rangle , \quad \forall \ \uom\, \widetilde{\varphi}(r,\uom) \in \mcW
$$
In $\mR^m$ we consider the space $\Omega(\mR^m) = \{ \uom\, \varphi(\ux) : \varphi(\ux) \in \mcD(\mR^m)  \}$. Clearly the functions in $\Omega(\mR^m)$ are no longer differentiable in the whole of $\mR^m$, since they are not defined at the origin due to the function $\uom = \frac{\ux}{|\ux|}$. By definition there is a one--to--one correspondence between the spaces $\mcD(\mR^m)$ and $\Omega(\mR^m)$.\\
For each $\widetilde{U}(r,\uom) \in \mcW'$ we define $U(\ux)$ by
$$
\langle \ U(\ux) , \uom\, \varphi(\ux) \ \rangle =  \langle \ \widetilde{U}(r,\uom) , \uom\, \widetilde{\varphi}(r,\uom)) \ \rangle , \quad \forall \ \uom\,\varphi(\ux) \in \Omega(\mR^m)
$$
Clearly $U(\ux)$ is a bounded linear functional on $\Omega(\mR^m)$, which we call a {\em signumdistribution}.\\

\noindent
Now start with a standard distribution $T(\ux) \in \mcD'(\mR^m)$ and let $\mT(r,\uom) \in \mcD'(\mR \times \mS^{m-1})$ be one of its spherical representations. Put
$\mS(r,\uom) = \uom\, \mT(r,\uom)$ which in its turn leads to the signumdistribution $S(\ux) \in \Omega(\mR^m)$. Then we consecutively have
\begin{eqnarray*}
\langle \  S(\ux) ,  \uom\, \varphi(\ux)  \ \rangle = \langle \  \mS(r,\uom) , \uom\, \widetilde{\varphi}(r,\uom)    \ \rangle &=& \langle \  \uom\, \mT(r,\uom) ,  \uom\,  \widetilde{\varphi}(r,\uom)  \ \rangle\\ &=& -\,  \langle \   \mT(r,\uom) ,   \widetilde{\varphi}(r,\uom)  \ \rangle = -\, \langle \  T(\ux) ,  \varphi(\ux)  \ \rangle 
\end{eqnarray*}
since $\uom^2 = -1$, and we call $S(\ux)$ a signumdistribution associated to the distribution $T(\ux)$ and denote it by $\uom\, T(\ux)$. It should be emphasized that for a given distribution $T(\ux)$ the associated signumdistribution $\uom\, T(\ux)$ is not uniquely defined but instead depends on the spherical representation of $T(\ux)$ chosen; moreover $\uom\, T$ is a mere notation, not the product of $T$ and $\uom$.

\begin{example}
{\em
A locally integrable function $f(\ux) \in L_1^{loc}(\mR^m)$ gives rise to a {\em regular distribution} $T_f$ via
$$
\langle \, T_f , \varphi(\ux) \ \rangle = \int_{\mR^m}\, f(\ux) \, \varphi(\ux)\, d\ux, \quad \forall\,  \varphi(\ux) \in \mcD(\mR^m)
$$
As $|\uom| = 1$ it is clear that also $f(\ux)\, \uom$ belongs to $ L_1^{loc}(\mR^m)$, thus generating the regular distribution $T_{f\uom}$ via
$$
\langle \, T_{f\uom} , \varphi(\ux) \ \rangle = \int_{\mR^m}\, f(\ux) \, \uom\,  \varphi(\ux)\, d\ux, \quad \forall\,  \varphi(\ux) \in \mcD(\mR^m)
$$
However the same integrals are defining the  {\em regular signumdistributions} $U_{f\uom}$ and $U_f$ by
$$
\langle \ U_{f\uom} , \uom\, \varphi(\ux) \ \rangle = -\, \int_{\mR^m}\, f(\ux) \, \varphi(\ux)\, d\ux, \quad \forall\,  \uom\, \varphi(\ux) \in \Omega(\mR^m)
$$
and
$$
\langle \, U_{f} , \uom\, \varphi(\ux) \ \rangle = \int_{\mR^m}\, f(\ux) \, \uom\,  \varphi(\ux)\, d\ux, \quad \forall\,  \uom\, \varphi(\ux) \in \Omega(\mR^m)
$$
A spherical representation of $T_f$ and $T_{f\uom}$ respectively is given by
$$
\langle \ \mT_f(r,\omega) , \Xi(r,\uom) \ \rangle = \int_{0}^\infty\, r^{m-1}\, dr\, \int_{\mS^{m-1}} f(r\uom)\, \Xi(r,\uom)\, dS_{\uom}
$$
$$
\langle \ \mT_{f\uom}(r,\omega) , \Xi(r,\uom) \ \rangle = \int_{0}^\infty\, r^{m-1}\, dr\, \int_{\mS^{m-1}} f(r\uom)\, \uom\, \Xi(r,\uom)\, dS_{\uom}
$$
since restricting to the space $\mcV$ leads to
$$
\langle \ \mT_f(r,\omega) , \widetilde{\varphi}(r,\uom) \ \rangle = \int_{0}^\infty\, r^{m-1}\, dr\, \int_{\mS^{m-1}} f(r\uom)\, \varphi(r\uom)\, dS_{\uom} = \langle \, T_f , \varphi(\ux) \ \rangle
$$
and
$$
\langle \ \mT_{f\uom}(r,\omega) , \widetilde{\varphi}(r,\uom) \ \rangle = \int_{0}^\infty\, r^{m-1}\, dr\, \int_{\mS^{m-1}} f(r\uom)\, \uom\, \varphi(r\uom)\, dS_{\uom} = \langle \, T_{f\uom} , \varphi(\ux) \ \rangle
$$
These particular spherical representations induce signumdistributions  associated to $T_f(\ux)$ and $T_{f\uom}(\ux)$, which we define to be $\uom\, T_f(\ux)$ and $\uom\, T_{f\uom}(\ux)$ respectively. It thus holds that
$$
\langle \  \uom\, T_f(\ux) ,  \uom\, \varphi(\ux)  \ \rangle =  -\, \langle \  T_f(\ux) ,  \varphi(\ux)  \ \rangle = -\,  \int_{\mR^m}\, f(\ux) \, \varphi(\ux)\, d\ux
$$
and
$$
\langle \  \uom\, T_{f\uom}(\ux) ,  \uom\, \varphi(\ux)  \ \rangle =  -\, \langle \  T_{f\uom}(\ux) ,  \varphi(\ux)  \ \rangle = -\,  \int_{\mR^m}\, f(\ux) \, \uom\, \varphi(\ux)\, d\ux
$$
Clearly $\uom\, T_f = U_{f\uom}$ and $\uom\, T_{f\uom} = -\, U_f$.\\
}
\end{example}

\begin{example}
{\em
Consider the distribution $T(\ux) =  \delta(\ux)$. Our aim is to define the signumdistributions $\uom\, \delta(\ux)$, $\p_r\, \delta(\ux)$ and $r\, \delta(\ux)$.\\
\noindent
A spherical representation of the delta distribution is given by
$$
\langle \ \mT(r,\omega) , \Xi(r,\uom) \ \rangle = \Sigma^0[\Xi(r,\uom)]\}|_{r=0}
$$
Indeed, when restricting to the space $\mcV$ and taking into account property (\ref{sigmanul}) we obtain
$$
\langle \ \mT(r,\omega) , \widetilde{\varphi}(r,\uom) \ \rangle =  \Sigma^0[\varphi(r\,\uom)]\}|_{r=0} = \langle \ \delta(\ux) , \varphi(\ux) \ \rangle
$$
This particular spherical representation of $T(\ux)$ induces a signumdistribution associated to $\delta(\ux) $, which we define to be  $\uom\, \delta(\ux)$ . It thus holds that
\begin{equation}
\label{omegadelta}
\langle \  \uom\, \delta(\ux) ,  \uom\, \varphi(\ux)  \ \rangle =  -\, \langle \  \delta(\ux) ,  \varphi(\ux)  \ \rangle 
\end{equation}

\noindent
Now consider the distribution $T_1(\ux) = (\uom\, \p_r)\, \delta(\ux) = \dirac\, \delta(\ux)$ (see Proposition \ref{omegadr}). A spherical representation of this distribution $T_1(\ux)$ is given by
$$
\langle \ \mT_1(r,\omega) , \Xi(r,\uom) \ \rangle = -\, m\, \{\p_r\, \Sigma^1[\Xi(r,\uom)]\}|_{r=0}
$$
Indeed, when restricting to the space $\mcV$ and taking into account property (\ref{sigmaeen}) we obtain
$$
\langle \ \mT_1(r,\omega) , \widetilde{\varphi}(r,\uom) \ \rangle =  -\, m\, \{\p_r\, \Sigma^1[\varphi(r\,\uom)]\}|_{r=0} = \langle \ \dirac\, \delta(\ux) , \varphi(\ux) \ \rangle
$$
This particular spherical representation of $T_1(\ux)$ induces a signumdistribution  associated to $T_1(\ux) = (\uom\, \p_r)\, \delta(\ux) $, which we define to be $-\, \p_r\, \delta(\ux)$. It thus holds that
\begin{equation}
\label{drdelta}
\langle \   \p_r\, \delta(\ux) ,  \uom\, \varphi(\ux)  \ \rangle =   \langle \  (\uom\, \p_r)\, \delta(\ux) ,  \varphi(\ux)  \ \rangle 
\end{equation}

\noindent
Finally as $\ux\, \delta(\ux) = 0$, we define the signumdistribution $r\, \delta(\ux)$ to be the zero signumdistribution.
}
\end{example}

\begin{example}
{\em
We consider the distribution $\dirac\, \delta(\ux) = (\uom\, \p_r)\, \delta(\ux)$.  As $\ux\, \dirac\, \delta(\ux) = m\, \delta(\ux)$, we first define, on the basis of a similar reasoning as in the previous example, the signumdistribution $r\, \dirac\, \delta(\ux)$ by
$$
\langle \, r\, \dirac\, \delta(\ux) , \uom\, \varphi(\ux) \ \rangle = m\, \langle \delta(\ux) , \varphi(\ux) \ \rangle
$$
In view of (\ref{omegadelta}) it clearly holds that $r\, \dirac\, \delta(\ux) = (-m)\, \uom\, \delta(\ux)$ or 
\begin{equation}
\label{romegadr}
r\ (\uom\ \p_r)\, \delta(\ux) = (-m)\, \uom\, \delta(\ux)
\end{equation}

\noindent
More generally, by considering the distribution
$$
\ux\, \dirac^{2k+1}\, \delta(\ux) = (m+2k)\, \dirac^{2k}\, \delta(\ux)
$$
we define the signumdistribution $r\, \dirac^{2k+1}\, \delta(\ux)$ by
$$
\langle \ r\, \dirac^{2k+1}\, \delta(\ux) , \uom\, \varphi(\ux) \ \rangle = (m+2k)\,  \langle \ \dirac^{2k}\, \delta(\ux) , \varphi(\ux) \ \rangle
$$
or, invoking the formulae obtained in Corollary \ref{omegaformulae}, 
$$
\langle \ r\, (\uom\, \p_r^{2k+1})\, \delta(\ux) , \uom\, \varphi(\ux) \ \rangle = (m+2k)\,  \langle \ \p_r^{2k}\, \delta(\ux) , \varphi(\ux) \ \rangle
$$

\noindent
Now we define the signumdistribution $\uom\, \dirac\, \delta(\ux) = \uom\, (\uom\, \p_r)\, \delta(\ux)$ to be given by
$$
\langle \ \uom\, \dirac\, \delta(\ux) , \uom\, \varphi(\ux) \ \rangle = -\, \langle \ \dirac\, \delta(\ux) , \varphi(\ux) \ \rangle
$$
In view of (\ref{drdelta}) it clearly holds that
$$
\uom\, (\uom\, \p_r)\, \delta(\ux) = - \, \p_r\, \delta(\ux)
$$

\noindent
More generally we define the signum distribution $\uom\, (\uom\, \p_r^{2k+1})\, \delta(\ux)$ by
$$
\langle \ \uom\, (\uom\, \p_r^{2k+1})\, \delta(\ux) , \uom\, \varphi(\ux) \ \rangle = -\, \langle \ (\uom\, \p_r^{2k+1})\, \delta(\ux) , \varphi(\ux) \ \rangle
$$

\noindent
Finally we define the signumdistribution $\p_r\, (\uom\, \p_r^{2k+1})\, \delta(\ux)$ by
$$
\langle \ \p_r\, (\uom\, \p_r^{2k+1})\, \delta(\ux) , \uom\, \varphi(\ux) \ \rangle = -\, \langle \ \p_r^{2k+2}\, \delta(\ux) , \varphi(\ux) \ \rangle
$$

}
\end{example}

\begin{example}
{\em
We consider the distribution $\Delta\, \delta(\ux) = -\, \dirac^2 \delta(\ux)$.
As $\ux\, \Delta\, \delta(\ux) = - 2\, \dirac\, \delta(\ux)$, we first define the signumdistribution $r\, \Delta(\ux)$ by
$$
\langle \ r\, \Delta\, \delta(\ux) , \uom\, \varphi(\ux) \ \rangle = - 2\, \langle \ \dirac \delta(\ux) , \varphi(\ux) \ \rangle
$$
Clearly $r\, \Delta\, \delta(\ux) = -2\, \p_r\, \delta(\ux)$, from which it also follows, by means of the results in Proposition \ref{partslaplace}, that  
$$
r\left(\invr\, \p_r\right) \delta(\ux) = \p_r\, \delta(\ux)
$$
and that
\begin{equation}
\label{rdrtwee}
r\, (\p_r^2)\, \delta(\ux) = -(m+1)\, \p_r\, \delta(\ux)
\end{equation}

\noindent
More generally, based upon the formula $\ux\, \dirac^{2k}\, \delta(\ux) = (2k)\, \dirac^{2k-1}\, \delta(\ux)$, we define the signumdistribution $r\, \dirac^{2k}\, \delta(\ux)$ by
$$
\langle \ r\, \dirac^{2k}\, \delta(\ux) , \uom\, \varphi(\ux) \ \rangle = (2k)\, \langle \ \dirac^{2k-1}\, \delta(\ux) , \varphi(\ux) \ \rangle
$$
or, again invoking the formulae obtained in Corollary \ref{omegaformulae},
$$
\langle \  r\, \p_r^{2k}\, \delta(\ux) , \uom\, \varphi(\ux) \ \rangle = -(m+2k-1)\, \langle \ (\uom\, \p_r^{2k-1})\, \delta(\ux) , \varphi(\ux) \ \rangle
$$

\noindent
Now we define the signumdistribution $\uom\, \dirac^{2k} \delta(\ux)$ by
$$
\langle \ \uom\, \dirac^{2k} \delta(\ux) , \uom\, \varphi(\ux) \ \rangle = -\, \langle \ \dirac^{2k} \delta(\ux) , \varphi(\ux) \ \rangle
$$
or
$$
\langle \ \uom\, \p_r^{2k} \delta(\ux) , \uom\, \varphi(\ux) \ \rangle = -\, \langle \ \p_r^{2k} \delta(\ux) , \varphi(\ux) \ \rangle
$$

\noindent
Finally we define the signumdistribution $\p_r\, \dirac^{2k} \delta(\ux)$ by
$$
\langle \ \p_r\, \dirac^{2k} \delta(\ux) , \uom\, \varphi(\ux) \ \rangle = \langle \ (\uom\, \p_r)\, \dirac^{2k} \delta(\ux) , \varphi(\ux) \ \rangle 
$$
which turns into
$$
\langle \ \p_r\, \p_r^{2k} \delta(\ux) , \uom\, \varphi(\ux) \ \rangle = \langle \ (\uom\, \p_r)\, \p_r^{2k} \delta(\ux) , \varphi(\ux) \ \rangle 
$$
or still
$$
\langle \ \p_r^{2k+1} \delta(\ux) , \uom\, \varphi(\ux) \ \rangle = \langle \ (\uom\, \p_r^{2k+1})\, \delta(\ux) , \varphi(\ux) \ \rangle 
$$
}
\end{example}

\begin{remark}
{\rm
For $k=0$ we obtain in particular that $r\, \delta	(\ux) = 0$. In fact this product is defined within the framework of standard distributions since the delta distribution is of finite order zero and the function $r$ is continuous in $\mR^m$.
}
\end{remark}

\begin{example}
{\em
Division of a standard distribution by a smooth function being allowed, we have
$$
\frac{1}{\ux}\, \delta(\ux) = \frac{1}{m}\, \dirac\, \delta(\ux) + \underline{S}_0
$$
where $\underline{S}_0$  stands for $\delta(\ux)\, \underline{c}_0$ with $\underline{c}_0$ an arbitrary constant vector, since $\ux\, \dirac\, \delta(\ux) = m\, \delta(\ux)$ and $\ux\, \underline{S}_0 = 0$. Making use of this formula we define the signumdistribution $\invr\, \delta(\ux)$ by
$$
\langle \ \invr\, \delta(\ux) , \uom\, \varphi(\ux) \ \rangle = \langle \ -\, \frac{1}{\ux}\, \delta(\ux) ,  \varphi(\ux) \ \rangle = \langle \ -\, \frac{1}{m}\, \dirac\, \delta(\ux) - \underline{S}_0 ,  \varphi(\ux) \ \rangle
$$
In view of
$$
\langle \ \p_r\, \delta(\ux) , \uom\, \varphi(\ux) \ \rangle = \langle \ \dirac\, \delta(\ux) , \varphi(\ux) \ \rangle
$$
and
$$
\langle \ \uom\, \underline{S}_0 , \uom\, \varphi(\ux) \ \rangle = \langle \ -\, \underline{S}_0 , \varphi(\ux) \ \rangle
$$
we obtain the following relation involving signumdistributions:
$$
\invr\, \delta(\ux)  = -\, \frac{1}{m}\, \p_r\, \delta(\ux) + \uom\, \delta(\ux)\, \underline{c}_0
$$
But as we expect the signumdistributions $\invr\, \delta(\ux)$ and  $\p_r\, \delta(\ux)$ to be SO$(m)$--invariant, the arbitrary vector constant $\underline{c}_0$ should be zero and we end up with
$$
\invr\, \delta(\ux)  = -\, \frac{1}{m}\, \p_r\, \delta(\ux)
$$
More generally we put
\begin{eqnarray*}
\langle \ \frac{1}{r^{2k+1}}\, \delta(\ux) , \uom\, \varphi(\ux) \, \rangle &=& \langle \ (-1)^{k+1}\, \frac{1}{\ux^{2k+1}}\, \delta(\ux) , \varphi(\ux) \ \rangle\\
&=& (-1)^{k+1}\, \langle \  \frac{1}{2^k k! (m+2k)(m+2k-2)\cdots(m)}\, \dirac^{2k+1}\delta(\ux) + \underline{S}_{2k}, \varphi(\ux) \ \rangle
\end{eqnarray*}
with $\underline{S}_{2k}$ an arbitrary vector linear combination of derivatives of the delta distribution up to order $2k$, which eventually leads  to
$$
\frac{1}{r^{2k+1}}\, \delta(\ux) = - \frac{(m-1)!}{(m+2k)!}\, \p_r^{2k+1}\, \delta(\ux) 
$$
Also, based upon the following formulae for division in cartesian coordinates:
\begin{eqnarray*}
\frac{1}{\ux}\, \dirac^{2k+1} \delta(\ux) &=& \frac{1}{2k+2}\, \dirac^{2k+2} \delta(\ux) + S_0\\
\frac{1}{\ux}\, \dirac^{2k} \delta(\ux) &=& \frac{1}{m+2k}\, \dirac^{2k+1} \delta(\ux) + \underline{S}_0
\end{eqnarray*}
and in a similar way as above, it is shown that
$$
\invr\, (\uom\, \p_r^{2k+1}) \delta(\ux) = - \frac{1}{m+2k+1}\, (\uom\, \p_r^{2k+2}) \delta(\ux) + \uom\ \delta(\ux) c_0
$$
and
$$
\invr\, \p_r^{2k} \delta(\ux) = - \frac{1}{m+2k}\, \p_r^{2k+1} \delta(\ux)
$$
}
\end{example}


\section{Some calculus}


There is a fundamental sequence of derivatives of the delta distribution, which are alternatively scalar and vector valued, generated by the action of the operator $(\uom\, \p_r)$:
$$
\delta \rightarrow (\uom\, \p_r)\,\delta \rightarrow  \p_r^2\,\delta \rightarrow  \ldots \rightarrow  (\uom\, \p_r^{2k-1})\,\delta \rightarrow  \p_r^{2k}\,\delta  \rightarrow  (\uom\, \p_r^{2k+1})\,\delta  \rightarrow \ldots
$$
and for each of the distributions in this sequence we have, in the examples of the foregoing section, defined, through the action of $\uom$, a specific associated signumdistribution, yielding in this way a parallel sequence of signumdistributions:
$$
\uom\, \delta \rightarrow  \p_r\,\delta \rightarrow  (\uom\, \p_r^2)\,\delta \rightarrow \ldots \rightarrow  \p_r^{2k-1}\,\delta \rightarrow  (\uom\, \p_r^{2k})\,\delta  \rightarrow  \p_r^{2k+1}\,\delta  \rightarrow \ldots
$$
Let us recall these definitions. The initial definition is the following.

\begin{definition}
\begin{eqnarray}
\label{def11}
\langle \ \uom\, \p_r^{2k}\, \delta(\ux) , \uom\, \varphi(\ux) \ \rangle &=& \langle \ -\, \p_r^{2k}\, \delta(\ux) , \varphi(\ux) \ \rangle\\[2mm]
\label{def12}
\langle \ \uom\, (\uom\, \p_r^{2k-1}\, \delta(\ux)) , \uom\, \varphi(\ux) \ \rangle &=& \langle \ -\, (\uom\, \p_r^{2k-1})\, \delta(\ux) , \varphi(\ux) \ \rangle
\end{eqnarray}
\end{definition}

\noindent
Whereupon we introduce the signumdistribution $\p_r^{2k-1}\, \delta(\ux)$ by

\begin{definition}
$$
\p_r^{2k-1}\, \delta(\ux) = -\, \uom\, (\uom\, \p_r^{2k-1}\, \delta(\ux))
$$
\end{definition}

\noindent
such that definition (\ref{def12}) may be rephrased as
\begin{equation}
\label{def12dash}
\langle \ \p_r^{2k-1}\, \delta(\ux) , \uom\, \varphi(\ux) \ \rangle = \langle \  (\uom\, \p_r^{2k-1})\, \delta(\ux) , \varphi(\ux) \ \rangle
\end{equation}

\noindent
There are two other actions on each of the distributions of the first sequence yielding a signumdistribution of the second sequence, viz. the actions by $r$ and by $\p_r$. Indeed, by an appropriate combination of the above definitions, we obtain the following calculus rules.

\begin{property}
One has
\begin{eqnarray*}
r\, (\uom\, \p_r^{2k+1})\, \delta(\ux) &=& -(m+2k)\, \uom\, \p_r^{2k}\, \delta(\ux)\\[2mm]
r\, \p_r^{2k}\, \delta(\ux) &=& -(m+2k-1)\, \p_r^{2k-1}\, \delta(\ux)\\[2mm]
\p_r\, (\uom\, \p_r^{2k+1})\, \delta(\ux) &=&  \uom\, \p_r^{2k+2}\, \delta(\ux)\\[2mm]
\p_r\, \p_r^{2k}\, \delta(\ux) &=&  \p_r^{2k+1}\, \delta(\ux)
\end{eqnarray*}
\end{property}

\noindent
One may wonder if there are actions transforming the signumdistributions from the second sequence back into distributions from the first sequence and the answer is positive. Indeed, the same actions apply on the signumdistributions from the second sequence. The basic action is again through the operator $\uom$, which yields the following definitions.

\begin{definition}
\begin{eqnarray}
\label{def21}
\langle \ \uom\, (\uom\, \p_r^{2k})\, \delta(\ux) , \varphi(\ux) \ \rangle &=& \langle \ \uom\, \p_r^{2k}\, \delta(\ux) , \uom\, \varphi(\ux) \ \rangle\\[2mm]
\label{def22}
\langle \ \uom\,  (\p_r^{2k-1}\, \delta(\ux)) , \varphi(\ux) \ \rangle &=& \langle \  \p_r^{2k-1}\, \delta(\ux) , \uom\, \varphi(\ux) \ \rangle
\end{eqnarray}
\end{definition}

\noindent
Comparing the definitions (\ref{def21}) and (\ref{def11}) it is clear that the distribution $\uom\, (\uom\, \p_r^{2k})\, \delta(\ux)$ is nothing else but the distribution 
$ -\, \p_r^{2k}\, \delta(\ux)$, while comparing definitions (\ref{def22}) and (\ref{def12dash}) shows that the distribution $\uom\,  (\p_r^{2k-1}\, \delta(\ux))$ is indeed the distribution $(\uom\, \p_r^{2k-1})\, \delta(\ux)$.\\

\noindent
For the actions of the operators $r$ and $\p_r$ on the signumdistributions, which are defined in a similar way as the actions of $r$ and $\p_r$ on distributions, we obtain the following computation rules.

\begin{property}
One has
\begin{eqnarray*}
r\, \p_r^{2k+1}\, \delta(\ux) &=& -(m+2k)\, \p_r^{2k}\, \delta(\ux)\\[2mm]
r\, (\uom\, \p_r^{2k})\, \delta(\ux) &=& -(m+2k-1)\, (\uom\, \p_r^{2k-1})\, \delta(\ux)\\[2mm]
r\, \left(\invr\, \p_r^{2k} \delta(\ux)\right) &=& \p_r^{2k} \delta(\ux)\\[1mm]
r\, \left(\invr\, \uom\, \p_r^{2k+1} \delta(\ux)\right) &=& \uom\, \p_r^{2k+1} \delta(\ux)\\[2mm]
\p_r\, (\p_r^{2k+1})\, \delta(\ux) &=&  \p_r^{2k+2}\, \delta(\ux)\\[2mm]
\p_r\, (\uom\, \p_r^{2k})\, \delta(\ux) &=&  (\uom\, \p_r^{2k+1})\, \delta(\ux)
\end{eqnarray*}
\end{property}

This leads to the following completely symmetric picture
$$
\begin{array}{ccccccccc}
  \ldots & \longrightarrow & (\uom\, \p_r^{2k-1})\,\delta & \longrightarrow  & \p_r^{2k}\,\delta & \longrightarrow & (\uom\, \p_r^{2k+1})\,\delta & \longrightarrow & \ldots \\[1mm]
  &&\hspace{-4.5mm}{}_{\uom}&\hspace{2mm}{}_r \hspace{7mm}{}_{\p_r}&\hspace{-4.5mm}{}_{\uom}&\hspace{2mm}{}_r \hspace{7mm}{}_{\p_r}&\hspace{-4.5mm}{}_{\uom}&&\\[-2mm]
  && \uparrow &  \nwarrow \hspace{-0.4mm} \nearrow & \uparrow &  \nwarrow \hspace{-0.4mm} \nearrow & \uparrow &&\\[-1.1mm]
   && \downarrow & \swarrow \hspace{-0.4mm} \searrow & \downarrow & \swarrow \hspace{-0.4mm} \searrow & \downarrow &&\\[-2mm]
&&\hspace{-4.5mm}{}^{\uom}&\hspace{2mm}{}^r \hspace{7mm}{}^{\p_r}&\hspace{-4.5mm}{}^{\uom}&\hspace{2mm}{}^r \hspace{7mm}{}^{\p_r}&\hspace{-4.5mm}{}^{\uom}&&\\[1mm]
 \ldots & \longrightarrow  & \p_r^{2k-1}\,\delta & \longrightarrow & (\uom\, \p_r^{2k})\,\delta & \longrightarrow  & \p_r^{2k+1}\,\delta & \longrightarrow & \ldots
 \end{array}
$$

\begin{remark}
{\rm
When composing two operators out of $r$, $\p_r$ and $\uom$, six operators originate: $r^2$, $r\, \p_r$, $r\, \uom$, $\p_r^2$, $\uom\, \p_r$ and $\uom^2$, which are traditional operators whose actions on distributions are well--defined. This means that the consecutive action by any two of the operators $r$, $\p_r$ and $\uom$ should lead to a known result, which is a serious test for all calculus rules established above. We now prove that this is indeed the case.\\

\noindent
(i) By the calculus rules we have
$$
r^2 (\uom\, \p_r^{2k+1}\, \delta) = - (m+2k) r (\uom\, \p_r^{2k}\, \delta) = (m+2k-1)(m+2k) (\uom\, \p_r^{2k-1})\, \delta
$$
and
$$
r^2 (\p_r^{2k}\, \delta) = - (m+2k-1) r (\p_r^{2k-1}\, \delta) = (m+2k-2)(m+2k-1) \p_r^{2k-2}\, \delta
$$
On the other hand, invoking the identities
\begin{eqnarray}
\label{identityodd}
\ux^2 \dirac^{2k+1}\, \delta(\ux) &=& (m+2k)(2k)\, \dirac^{2k-1}\, \delta(\ux)\\
\label{identityeven}
\ux^2 \dirac^{2k}\, \delta(\ux) &=& (m+2k-2)(2k)\, \dirac^{2k-2}\, \delta(\ux)
\end{eqnarray}
and the formulae of Corollary \ref{omegaformulae}
we have
\begin{eqnarray*}
r^2 (\uom\, \p_r^{2k+1}\, \delta) &=& -  \frac{(-1)^k}{2^k k!}\, (m+1)(m+3)\cdots(m+2k-1)\, \ux^2 \dirac ^{2k+1}\, \delta\\
&=& (m+2k-1)(m+2k) (\uom\, \p_r^{2k-1})\, \delta
\end{eqnarray*}
and
\begin{eqnarray*}
r^2 (\p_r^{2k}\, \delta) &=& -  \frac{(-1)^k}{2^k k!}\, (m+1)(m+3)\cdots(m+2k-1)\, \ux^2 \dirac ^{2k}\, \delta\\
&=& (m+2k-2)(m+2k-1) (\uom\, \p_r^{2k-2})\, \delta
\end{eqnarray*}

\noindent
(ii) The Euler operator measures the degree of homogeneity and thus
\begin{eqnarray*}
r\, \p_r\, (\uom\, \p_r^{2k+1}\, \delta) &=& -(m+2k+1)\, \uom\, \p_r^{2k+1}\, \delta\\
r\, \p_r\, (\p_r^{2k}\, \delta) &=& -(m+2k) \p_r^{2k}\, \delta
\end{eqnarray*}
while the calculus rules lead to
\begin{eqnarray*}
r\, \p_r\, (\uom\, \p_r^{2k+1}\, \delta) &=& r (\uom\, \p_r^{2k+2}\, \delta) =  -(m+2k+1)\, \uom\, \p_r^{2k+1}\, \delta\\
r\, \p_r\, (\p_r^{2k}\, \delta) &=& r (\p_r^{2k+1}\, \delta) = -(m+2k) \p_r^{2k}\, \delta
\end{eqnarray*}

\noindent
(iii) By the calculus rules we obtain
\begin{eqnarray*}
r\, \uom\, (\uom\, \p_r^{2k+1}\, \delta) &=&  -\, r\, \p_r^{2k+1}\, \delta = (m+2k)\, \p_r^{2k}\, \delta \\
r\, \uom\, (\p_r^{2k}\, \delta) &=& - (m+2k-1)\, \uom\, \p_r^{2k-1}\, \delta
\end{eqnarray*}
while invoking the identities (\ref{identityodd}) and (\ref{identityeven}) respectively leads to
\begin{eqnarray*}
r\, \uom\, (\uom\, \p_r^{2k+1}\, \delta) &=&    \frac{(-1)^k}{2^k k!}\,   (m+1)(m+3)\cdots(m+2k-1)\, \ux\,  \dirac^{2k+1}\, \delta     \\
&=& (m+2k)\, \p_r^{2k}\, \delta 
\end{eqnarray*}
and
\begin{eqnarray*}
r\, \uom\, (\p_r^{2k}\, \delta) &=&    \frac{(-1)^k}{2^k k!}\,   (m+1)(m+3)\cdots(m+2k-1)\, \ux\,  \dirac^{2k}\, \delta     \\
&=& -(m+2k-1)\, \uom\, \p_r^{2k-1}\, \delta 
\end{eqnarray*}

\noindent
(iv) The calculus rules lead to
\begin{eqnarray*}
\p_r^2\, (\uom\, \p_r^{2k+1}\, \delta) &=& \uom\, \p_r^{2k+3}\, \delta\\
\p_r^2\, (\p_r^{2k}\, \delta) &=& \p_r^{2k+2}\, \delta
\end{eqnarray*}
On the other hand we can make use of the identities
\begin{eqnarray*}
(\uom\, \p_r)\, \dirac^{2k}\, \delta(\ux) &=& \dirac^{2k+1}\, \delta(\ux) \\
(\uom\, \p_r)\, \dirac^{2k+1}\, \delta(\ux) &=& \frac{m+2k+1}{2(k+1)}\, \dirac^{2k+2}\, \delta(\ux)
\end{eqnarray*}
to obtain
\begin{eqnarray*}
\p_r^2\, (\uom\, \p_r^{2k+1}\, \delta) &=&  \frac{(-1)^k}{2^k k!}\, (m+1)(m+3)\cdots(m+2k-1)\, (-1)\, \frac{m+2k+1}{2(k+1)}\, \dirac^{2k+3}\, \delta\\
&=& \frac{(-1)^{k+1}}{2^{k+1} (k+1)!}\, (m+1)(m+3)\cdots(m+2k-1)(m+2k+1)\, \dirac^{2k+3}\, \delta = \uom\, \p_r^{2k+3}\, \delta
\end{eqnarray*}
and
\begin{eqnarray*}
\p_r^2\, (\p_r^{2k}\, \delta) &=&  \frac{(-1)^k}{2^k k!}\, (m+1)(m+3)\cdots(m+2k-1)\, (-1)\, \frac{m+2k+1}{2(k+1)}\, \dirac^{2k+2}\, \delta\\
&=& \frac{(-1)^{k+1}}{2^{k+1} (k+1)!}\, (m+1)(m+3)\cdots(m+2k-1)(m+2k+1)\, \dirac^{2k+2}\, \delta =  \p_r^{2k+2}\, \delta
\end{eqnarray*}

\noindent
(v) On the one hand we have by the calculus rules
\begin{eqnarray*}
\uom\, \p_r\, (\uom\, \p_r^{2k+1}\, \delta) &=& \uom\, (\uom\, \p_r^{2k+2}\, \delta) = -\, \p_r^{2k+2}\, \delta\\
\uom\, \p_r\, (\p_r^{2k}\, \delta) &=& \uom\, (\p_r^{2k+1}\, \delta) = \uom\, \p_r^{2k+1}\, \delta
\end{eqnarray*}
and on the other
\begin{eqnarray*}
\uom\, \p_r\, (\uom\, \p_r^{2k+1}\, \delta) &=&  \uom\, \p_r\, (-1)^k (\uom\, \p_r)^{2k+1}\, \delta =   (-1)^k (\uom\, \p_r)^{2k+2}\, \delta = -\, \p_r^{2k+2}\, \delta\\
\uom\, \p_r\, (\p_r^{2k}\, \delta) &=& \uom\, \p_r\, (-1)^k (\uom\, \p_r)^{2k}\, \delta = (-1)^k (\uom\, \p_r)^{2k+1}\, \delta = \uom\, \p_r^{2k+1}\ \delta
\end{eqnarray*}

\noindent
(vi) The action by $\uom^2 = -1$ is trivial.

}
\end{remark}


\section{Conclusion}


In our quest for an unambiguous meaningful definition, in the framework of spherical coordinates, of the radial derivative $\p_r\, \delta(\ux)$ of the delta distribution in $\mR^m$, and of other ``forbidden'' actions on the delta distribution such as $r\, \delta(\ux)$ and $\uom\, \delta(\ux)$, we were faced with the impossibility to achieving this within the familiar setting of the traditional distributions. Instead we had to introduce a new space of continuous linear functionals on a space of test functions showing a singularity at the origin, for which we coined the term {\em signumdistributions}, bearing in mind that $\uom = \frac{\ux}{|\ux|}$ may be interpreted as the higher dimensional counterpart to the {\em signum} function on the real line. It turns out that the actions by $r$, $\uom$ and $\p_r$ map a distribution to a signumdistribution and vice versa, and a number of efficient calculus rules for handling these transitions for the delta distribution were established. As the composition of any two operators from  $r$, $\uom$ and $\p_r$ results in a ``legal" and well--defined action on distributions, and on the delta distribution in particular, these calculus rules were positively tested for this phenomenon. Finally it should be mentioned that spaces of test functions showing a singularity at the origin also appear in the theory of so--called {\em thick distributions} (see \cite{estrada}). The possible relationships between the signumdistributions and the thick distributions are subject of current research.


\section{Acknowledgement}

The first author wants to thank Kevin Coulembier, Hendrik De Bie, Hennie De Schepper, and David Eelbode for their interest in and their valuable comments on
the topic treated in this paper.



\end{document}